\documentclass[11pt]{article}
\usepackage{amsmath,amsthm,amsfonts,latexsym,amssymb,enumerate,color}
\usepackage{graphicx}

\usepackage{color}

\newcommand\dist{\mathop{\rm dist}\nolimits}

\newtheorem{thm}{Theorem}[section]
\newtheorem{prop}[thm]{Proposition}

\newtheorem{cor}[thm]{Corollary}
   
\newtheorem{exam}[thm]{Example}
\newtheorem{defn}[thm]{Definition}

\newcommand\beq{\begin{equation}}
\newcommand\eeq{\end{equation}}

\newcommand\CC{{\mathbb C}}

\newcommand\DD{{\mathbb D}}
\newcommand\cH{{\mathcal H}}
\newcommand\cK{{\mathcal K}}

\newcommand\cM{{\mathcal M}}

\newcommand{\TT}{{\mathbb T}}
\newcommand{\ZZ}{{\mathbb Z}}

\newcommand\rank{\mathop{\rm rank}\nolimits}

\def\spam{\mathop{\rm \overline{span}}\nolimits}
\def\epsilon{\varepsilon}
\def\codim{\mathop{\rm codim}\nolimits}
\def\inter{\mathop{\rm int}\nolimits}

\def\beginpf{\begin{proof}}
\def\endpf{\end{proof}}

\title{Inner functions and operator theory}

\author{Isabelle Chalendar\\ \small{Universit\'e Lyon 1, INSA de Lyon, Ecole Centrale de Lyon}\\ \small{ CNRS, UMR 5208, Institut Camille Jordan}\\ \small{ 43 bld. du 11 novembre 1918, F-69622 Villeurbanne Cedex, France}\\ \small{\tt E-mail: chalendar@math.univ-lyon1.fr}\and Pamela Gorkin\\ \small{Department of Mathematics, Bucknell University, Lewisburg, PA 17837, U.S.A.} \\ \small{\tt E-mail: pgorkin@bucknell.edu} \and and \\ Jonathan R. Partington\\\small{School of Mathematics, University of Leeds, Leeds LS2 9JT, U.K.} \\ \small{\tt E-mail: J.R.Partington@leeds.ac.uk}}

\begin{document}

\maketitle

\begin{abstract}
This tutorial paper presents a survey of results, both classical and new, linking inner functions and operator theory.
Topics discussed include invariant subspaces, universal operators, Hankel and Toeplitz operators, model spaces,
truncated Toeplitz operators, restricted shifts, numerical ranges, and interpolation.
\end{abstract}

\section{Introduction}

Inner functions originally arose in the context of operator theory, via Beurling's theorem on the invariant subspaces of the
unilateral shift operator. Since then, they have been seen in numerous contexts in the theory of 
function spaces. This tutorial
paper surveys some of the many ways in which operators and inner functions are linked: these include the 
invariant subspace problem, the
theory of Hankel and Toeplitz operators and the rapidly-developing area of
model spaces and the operators acting on them. 

The paper is an expanded version of a mini-course given at
the
Eleventh Advanced Course in 
Operator Theory and Complex Analysis, held in 
Seville in June 2014.

\subsection{Hardy spaces and shift-invariant subspaces}

All our spaces will be complex. We write $\DD$ for the open
unit disc in $\CC$ and $\TT=\partial \DD$, the unit circle.

Recall that Hardy space $H^2$ or $H^2(\DD)$ is the space of analytic functions on   $\DD$ with square-summable Taylor coefficients; that is,
\[
H^2(\DD)=\{f: \DD \to \CC \hbox{ analytic}, f(z)=\sum_{n=0}^\infty a_nz^n, \|f\|^2= \sum_{n=0}^\infty |a_n|^2 < \infty \}.
\]
Also $H^2(\DD)$ embeds isometrically as a closed subspace of $L^2(\TT)$ via
\[
\sum_{n=0}^\infty a_n z^n \mapsto \sum_{n=0}^\infty a_n e^{int},
\]
where the series converges almost everywhere on $\TT$ as well as in the norm of $L^2(\TT)$. Indeed, $\lim_{r \to 1-} f(re^{it})$ exists almost
everywhere and gives the boundary values of a function $f$ in $H^2(\DD)$. (See, for
example \cite{hoffman}.)

It is useful to use the isometric isomorphism $\ell^2(\ZZ) \to L^2(\TT)$ given by $(a_n)_{n \in \ZZ} \mapsto \sum_{n=-\infty}^\infty a_n e^{int}$, which is a consequence of the
Riesz--Fischer theorem; this
restricts to an isomorphism $\ell_2(\ZZ_+) \to H^2(\DD)$.\\

The first connection between inner functions and operator theory arises on considering the right shift $R: \ell^2(\ZZ) \to \ell^2(\ZZ)$. We may ask what its closed invariant subspaces are;
that is, the subspaces $\cM \subset L^2(\TT)$ such that $R\cM\subset \cM$. The answer is to look at the unitarily equivalent operator $S$ of ``multiplication by $z$'' on $L^2(\TT)$.
\[
\begin{matrix}
\ell^2(\ZZ) & \stackrel{ R}{\rightarrow}& \ell^2(\ZZ) \\
\downarrow && \downarrow \\
L^2(\TT) & \stackrel{S}{\rightarrow} & L^2(\TT)
\end{matrix}
\]

There are two cases, for $\cM$ a nontrivial closed subspace of $L^2(\TT)$:\\

(i) $S\cM=\cM$, if and only if there is a measurable subset $E \subset \TT$ such that $\cM=\{f \in L^2(\TT): f_{|\TT \setminus E}=0 \hbox{ a.e.}\}$ (Wiener \cite[Ch.~II]{Wiener}).\\

(ii) $S\cM \subsetneq \cM$, if and only if there is a unimodular function $\phi \in L^\infty(\TT)$ such that $\cM=\phi H^2$ (Beurling--Helson \cite{helson}).\\

As a sketch proof of (ii), which will be the more important for us, take $\phi \in \cM \ominus S\cM$ with $\|\phi\|_2=1$. One can verify that $\phi$ is unimodular and
that $\cM=\phi H^2$.

\begin{cor}[Beurling's theorem {\cite{beurling}; see also \cite[Thm.~II.7.1]{Garnett}, \cite[Sec. A.1.3]{Nikolksi}}]
Let $\cM$ be a nontrivial closed subspace of $H^2$; then $S\cM \subset \cM$
if and only if $\cM=\theta H^2$ where $\theta$ is {\em inner}, that is $\theta \in H^2(\DD)$
with $|\theta(e^{it})|=1$ a.e.
\end{cor}
It is easily seen that $\theta$ is unique up to multiplication by
a  constant of modulus 1.\\

Now, any function $h \in H^2$, apart from the zero function, has a 
multiplicative factorization $h=\theta u$, where $\theta$ is inner, and $u$ is {\em outer\/}: 
Beurling showed that outer functions satisfy
\[
\spam\{u, Su, S^2u, S^3u, \ldots \}= H^2,
\]
and they therefore have
an operatorial interpretation, as cyclic vectors for the shift $S$.
The inner-outer factorization is unique up to multiplication by a constant  of modulus one. 

\subsection{Examples of inner functions}

If $\cM$ is a shift-invariant subspace of finite codimension, then $\theta$ is a 
finite Blaschke product,
\[
\theta(z) = \lambda    \prod_{j = 1}^n \frac{z-\alpha_j}{1 - \overline{\alpha_j}z},
\] 
with $|\lambda|=1$ and $\alpha_1,\ldots,\alpha_n \in \DD$. Then
\[
\cM= \{f \in H^2: f(\alpha_1)=\cdots=f(\alpha_n)=0\},
\]
with the obvious interpretation in the case of non-distinct $\alpha_j$.
We may also form infinite Blaschke products
\[
\theta(z) = \lambda z^p  \prod_{j = 1}^\infty \frac{|\alpha_j|}{\alpha_j}\frac{\alpha_j - z}{1 - \overline{\alpha_j}z},
\] 
where $|\lambda| = 1$, all the $\alpha_j$ lie in $\mathbb{D}\setminus\{0\}$, 
$p$ is a non-negative integer and $\sum_{j=1}^\infty (1-|\alpha_j|) < \infty$. Recall that the sequences of $\mathbb{D}$ satisfying the last condition are called 
Blaschke sequences.  

There is also a class of inner functions without zeroes, namely the
singular inner functions, which may be written as
\[
\theta(z)=\exp \left[ -\int_{-\pi}^{\pi} \frac{e^{it}+z}{e^{it}-z} \, d\mu(t) 
\right],
\]
where $\mu$ is a singular positive measure on $[-\pi,\pi)$.
For example if $\mu$ is a Dirac mass at $0$, then
$\theta(z)=\exp((z+1)/(z-1))$.

A complete description of inner functions is now available, as they are given as $Bs$, where
$B$ is a Blaschke product and $s$ is a singular inner function. Either factor may
be absent.

Note that if $\theta_1$ and $\theta_2$ are inner, then $\theta_1 \overline{\theta_2}$
is unimodular on $\TT$. These are not all the unimodular functions, but if
$\phi \in L^\infty(\TT)$ is unimodular then for each $\epsilon>0$ it can be
factorized as $\phi=h_1\overline{h_2}$, with $h_1,h_2 \in H^\infty$ 
and $\|h_1\|,\|h_2\| < 1+\epsilon$ (see \cite{bourgain,barclay,CPbook}).
Related to this is the Douglas--Rudin theorem that the quotients $ \theta_1 \overline{\theta_2}$
with $\theta_1$ and $\theta_2$ inner
are uniformly dense in the unimodular functions in $L^\infty(\TT)$ (see \cite{DR}).

Of particular importance are the interpolating Blaschke products: a Blaschke product $B$ with zeroes $(z_j)$ is interpolating if its zero sequence is an interpolating sequence for $H^\infty$ or, equivalently, there exists $\delta> 0$ such that 
$$\inf_{ k} \prod_{j: j \ne k} \left|\frac{z_j - z_k}{1 - \overline{z_k} z_j}\right| = \delta.$$ 
These Blaschke products play an important role in the study of bounded analytic functions: consider a closed subalgebra $B$ of $L^\infty$ containing $H^\infty$ properly. 
In establishing a conjecture of R. G. Douglas, Chang and Marshall \cite{c, m}  proved that such algebras (now called Douglas algebras) can be characterized using interpolating Blaschke products: 
if $$U_B = \{b: b~\mbox{interpolating and}~b^{-1} \in B\},$$ then an algebra is a Douglas algebra if and only if it is the closed algebra generated by $H^\infty$ and 
the conjugates of the functions in $U_B$. In other words,
$B = [H^\infty, \overline{U_B}]$. 
Much more is known about interpolating Blaschke products: in particular, P. Jones \cite{Jones} showed that one can take the Blaschke products in the Douglas--Rudin theorem to be interpolating. Related work can be found in \cite{ms}, \cite{gn}, and \cite{Garnett}. 
One very interesting question remains open: can every Blaschke product be approximated (uniformly) by an interpolating Blaschke product? Hjelle and Nicolau \cite{hn} have shown that given a Blaschke product, $B$, there is an interpolating Blaschke product that  approximates $B$ in modulus on $\mathbb{D}$, but this is the best result to date.

\section{Some operators associated with inner functions}

\subsection{Isometries}

(i) It is not hard to see that the {\em analytic Toeplitz operator}
or {\em Laurent operator}, $T_\phi: H^2 \to H^2$, $f \mapsto \phi f$, where $\phi \in H^\infty$,
is an isometry if and only if $\phi$ is inner. Moreover
$\codim \phi H^2 < \infty$ if and only if $\phi$ is a finite Blaschke product.\\

\noindent (ii) For $\phi: \DD \to \DD$ holomorphic, we may consider the 
{\em composition operator\/} $C_\phi: H^2 \to H^2$, $f \mapsto f \circ \phi$.
See for example \cite{cowen-mac} for full details on these. In particular,
by Littlewood's subordination theorem \cite{littlewood}, $C_\phi$ is automatically
continuous. 

It is a result of Nordgren \cite{nordgren} that $C_\phi$ is an isometry if and 
only if $\phi$ is inner and $\phi(0)=0$. Note that if $\phi$ is inner, with $\phi(0)=0$,
then for $n>m$ we have
\[
\langle \phi^n,\phi^m\rangle = 
\langle \phi^{n-m},1\rangle = \phi(0)^{n-m} = 0,
\]
so that the orthonormal sequence $(z^n)_{n \ge 0}$ in $H^2$ is mapped to the
orthonormal sequence $(\phi^n)_{n \ge 0}$.

Conversely, since $\langle z,1 \rangle=0$, we must have $\phi(0)=\langle \phi,1\rangle=0$ 
if $C_\phi$ is to be an isometry. Also the condition $\|\phi^n\|=1$ for all $n$ 
can be used to check that $\phi$ is inner.\\

Bayart \cite{bayart} shows that $C_\phi$ is {\em similar\/} to an isometry
if and only if $\phi$ is inner and $\phi(p)=p$ for some $p \in \DD$.

\subsection{Universal operators}

An operator $U$ defined on a separable infinite-dimensional Hilbert space $\cH$ is said to be {\em universal\/}
in the sense of Rota, if for every
operator $T$ on a Hilbert space $\cK$ there is a constant $\lambda \in \CC$ and an
invariant subspace $\cM$ for $U$ such that $T$ is similar to the restriction $\lambda U_{|\cM}$.

\[
\begin{matrix}
\cH & \stackrel{ \lambda U}{\rightarrow}& \cH \\
\uparrow && \uparrow \\
\cM & \stackrel{ \lambda U}{\rightarrow}& \cM\\
J\downarrow && J\downarrow \\
\cK & \stackrel{ T}{\rightarrow}& \cK\\
\end{matrix}
\]

The following theorem provides many examples of universal operators.

\begin{thm}[Caradus \cite{caradus}]\label{thm:caradus}
If the operator $H: \cH \to \cH$ is surjective with infinite-dimensional kernel,
then it is universal.
\end{thm}

(a) Take $\theta$ inner, but not a finite Blaschke product. Then using Theorem \ref{thm:caradus}
one can show that the
Toeplitz operator $T_{\overline\theta}=T^*_\theta: H^2 \to H^2$,
with $f \mapsto P_{H^2}(\overline\theta f)$ is universal.

Such an operator $T^*_\theta$ is similar to the backward shift $A$ on $L^2(0,\infty)$,
given by
\[
Af(t)=f(t+1),
\]
which by the Laplace transform is unitarily equivalent to the adjoint of the
operator $M_{e^{-s}}$ of multiplication of $e^{-s}$ on the Hardy space $H^2(\CC_+)$
of the right half-plane $\CC_+$ (here $s$ is the independent variable). Note that $e^{-s}$
is inner: still, in spite of Beurling's theorem mentioned above, there is no
usable characterization of the invariant subspaces of $A$.\\

(b) Let $\phi: \DD \to \DD$ be defined by
\[
\phi(z) = \frac{z+1/2}{1+z/2};
\]
this is a (hyperbolic) automorphism fixing $\pm 1$. The composition operator $C_\phi$ has spectrum
given by
\[
\sigma(C_\phi)= \{z \in \CC: 1/\sqrt3 \le |z| \le \sqrt 3\}.
\]
For $\lambda \in \inter \sigma(C_\phi)$, it can be shown that $C_\phi - \lambda I$ is
universal  \cite{NRW}. Note that it has the same invariant subspaces as $C_\phi$, and
a complete description of them would give a solution to the invariant subspace problem.

These ideas have stimulated studies on cyclic vectors and minimal invariant subspaces
for $C_\phi $ (e.g. \cite{mortini} and \cite{eva-pam}).

\subsection{Hankel and Toeplitz operators}

We begin with the orthogonal decomposition
\[
L^2(\TT) = H^2 \oplus \overline{H^2_0}
\]
into closed subspaces spanned by $\{e^{int}: n \ge 0\}$ and $\{e^{int}: n < 0\}$, respectively. Write $P:L^2(\TT) \to H^2$ for the
orthogonal projection.

\begin{defn}
Let $\phi \in L^\infty(\TT)$. Then the {\em Toeplitz operator\/} $T_\phi: H^2 \to H^2$ is defined by $T_\phi f = P(\phi f)$ for $f \in H^2$. The {\em Hankel operator\/}
$\Gamma_\phi: H^2 \to \overline{H^2_0}$ is defined by $\Gamma_\phi f=(I-P)\phi f$ for $f \in H^2$.
\end{defn}

It is well known that $\|T_\phi\|=\|\phi\|_\infty$ (see
Brown--Halmos \cite{BH}) and that $\|\Gamma_\phi\|=\dist(\phi,H^\infty)$ (see Nehari \cite{nehari}).

\subsection{Kernels}

(i) If $u \in \ker \Gamma_\phi$, then $\phi u \in H^2$, so that $z\phi u \in H^2$ and $zu \in \ker \Gamma_\phi$.
Hence, by Beurling's theorem, $\ker \Gamma_\phi = \theta H^2$ for some inner function $\theta$. 

For example, if $\theta$ is inner, then $u \in \ker \Gamma_{\overline\theta}$ if and only if $\overline\theta u \in H^2$,
which happens if and only if $u \in \theta H^2$. So all Beurling subspaces occur as Hankel kernels.\\

\noindent (ii) Suppose that $\theta$ is inner. Then
 $f \in \ker T_{\overline\theta}$ if and only if $\langle \overline\theta f,g\rangle=0$ for all $g \in H^2$. This is equivalent to the
condition $\langle f,\theta g\rangle =0$; that is,  $f \in H^2 \ominus \theta H^2$.
We shall study these spaces in Section \ref{sec:model}.

Toeplitz kernels in general
have the {\em near-invariance\/} property. If $u \in H^2$ and $\theta u \in \ker T_\phi$ for
some inner function $\theta$, then $\phi \theta u = \overline z\overline h$ for some $h \in H^2$.
Hence $\phi u = \overline\theta \overline z\overline h$ and thus $u \in \ker T_\phi$. 

That is, if $v \in \ker T_\phi$ and $v/\theta \in H^2$, then $v/\theta \in \ker T_\phi$.

In particular, if $v \in \ker T_\phi$ and $v/z \in H^2$, then $v/z \in \ker T_\phi$.
This property is not the same as being $S^*$-invariant, even though $S^*v=v/z$ if $v/z \in H^2$.

For example, let $\phi(z)=e^{-z}/z^2$. One may verify that
\[
\ker T_\phi=\{(a+bz)e^z: a,b \in \CC\}.
\]
 However $S^* e^z = \dfrac{e^z-1}{z}$,
which does not lie in $\ker T_\phi$.\\

Now Hitt \cite{hitt} showed that a subspace $\cM \subset H^2$ is nearly $S^*$-invariant if and only if it can be written as $\cM=fK_\theta$, where $\theta$ is inner, 
$\theta(0)=0$, $f \in \cM \ominus (\cM \cap zH^2)$,
and $K_\theta$ is the model space $H^2 \ominus \theta H^2$, discussed in Section \ref{sec:model}.\\

Moreover, Hayashi \cite{hayashi86,hayashi90} showed that such an $\cM$ is in fact a Toeplitz kernel
if and only if the function $f$ has the property that $f^2$ is {\em rigid}, which means that if $g \in H^1$ with $g/f^2>0$ a.e., then $g=\lambda f^2$ for
some constant $\lambda>0$. A rigid function is necessarily outer.

\section{Model spaces}
\label{sec:model}

\subsection{Definitions and examples}

Since the invariant subspaces for $S$ have the form $\theta H^2$, with $\theta$ inner, those for $S^*$
have the form $H^2 \ominus \theta H^2$, usually written $K_\theta$. Such spaces are
called {\em model spaces}.

\begin{exam}{\rm  (i) Take $\theta(z)=z^N$, which is inner. Then 
\[
K_\theta=\mathop{\rm span}\{1,z,z^2,\ldots,z^{N-1}\}.
\]

(ii) For $\theta(z)= \prod_{k=1}^N \dfrac{z-\alpha_k}{1-\overline{\alpha_k}z}$ with $\alpha_1,\ldots,\alpha_N$ distinct, we have
$f \in \theta H^2$ if and only if $f(\alpha_1)=\cdots=f(\alpha_N)=0$. Then 
\[
K_\theta=\mathop{\rm span}\left\{\frac{1}{1-\overline{\alpha_1}z},\ldots, \frac{1}{1-\overline{\alpha_N}z}\right\}.
\]
Indeed, for $\alpha \in \DD$, $k_\alpha: z \mapsto \dfrac{1}{1-\overline{\alpha}z}$ is the {\em reproducing kernel\/} at $\alpha$; i.e.,
\[
f(\alpha)= \langle f,k_\alpha \rangle \qquad \hbox{for} \quad f \in H^2,
\]
and clearly $f \in \theta H^2$ if and only if $f$ is orthogonal to $k_{\alpha_1},\ldots,k_{\alpha_N}$.\\

(iii) For a fixed $\tau>0$ we write 
\beq\label{eq:l2tau}
L^2(0,\infty)=L^2(0,\tau) \oplus L^2(\tau,\infty).
\eeq
 Under the Laplace transform this maps to the
orthogonal decomposition
\[
H^2(\CC_+) = K_\theta \oplus \theta H^2(\CC_+),
\]
where $\theta(s)=e^{-s\tau}$; that is, $\theta$ is inner. Then $K_\theta$ 
can be written as $e^{s\tau/2}PW_{\tau/2}$, where $PW_{\tau/2}$
is a {\em Paley--Wiener\/} space, 
consisting of entire functions, as considered 
in signal processing.
} 
\end{exam}

In general $K_\theta$ is finite-dimensional if and only if $\theta$ is a finite Blaschke product.\\

\subsection{Decompositions of $H^2$ and $K_B$.} Let $\theta$ be inner. Then
\[
H^2= K_\theta \oplus \theta K_\theta \oplus \theta^2 K_\theta \oplus \cdots.
\]
This is an orthogonal direct sum, since if $k_1, k_2 \in K_\theta$ and $0 \le m < n$, then
\[
\langle \theta^m k_1, \theta^n k_2 \rangle = \langle k_1, \theta^{n-m}k_2 \rangle =0,
\]
since $k_1 \perp \theta H^2$.

Note that $T_\theta$ acts as a shift here, i.e.,
\[
\theta(k_1+\theta k_2 + \theta^2 k_3 + \cdots)
=
\theta k_1+\theta^2 k_2 + \theta^3 k_3 + \cdots.
\]
A special case of this can be identified from (\ref{eq:l2tau}) above, since
\[
L^2(0,\infty)= L^2(0,\tau) \oplus L^2(\tau,2\tau) \oplus \cdots.
\]

We now look at model spaces corresponding to infinite Blaschke products. If $\alpha_1,\alpha_2,\ldots$ are the zeroes of an infinite Blaschke product $B$ (assumed distinct),
then an orthonormal basis of $K_B$ is the Takenaka--Malmquist--Walsh basis given by orthonormalizing the sequence of reproducing kernels associated with the $(\alpha_n)$.
We have
\begin{eqnarray*}
e_1(z) &=& \frac{(1-|\alpha_1|^2)^{1/2}}{1-\overline{\alpha_1}z},\\
e_2(z) &=& \frac{(1-|\alpha_2|^2)^{1/2}}{1-\overline{\alpha_2}z} \left(\frac{z-\alpha_1}{1-\overline{\alpha_1}z}\right),\\
\noalign{\hbox{and, in general}}\\
e_n(z) &=& \frac{(1-|\alpha_n|^2)^{1/2}}{1-\overline{\alpha_n}z}\left( \prod_{k=1}^{n-1}\frac{z-\alpha_k}{1-\overline{\alpha_k}z}\right).
\end{eqnarray*}
It is easily checked that these are orthonormal, and have the same closed span as the
reproducing kernels $\dfrac{1}{1-\overline{\alpha_1}z},\ldots,\dfrac{1}{1-\overline{\alpha_n}z},\ldots$. This
closed span  is $K_B$ when the $(\alpha_n)$ form a Blaschke sequence, and $H^2$
otherwise.

\subsection{Frostman's theorem and mappings between model spaces}

The following result shows that inner functions are not far from   Blaschke
products, in a precise sense.

\begin{thm}[Frostman \cite{frostman}]
Let $\theta$ be any inner function. Then, for $\alpha \in \DD$,
the function $\dfrac{\theta-\alpha}{1-\overline{\alpha}\theta}$ is also inner; it is a Blaschke product with distinct zeroes 
for all $\alpha \in \DD$ outside an exceptional set
$E$ such that for each $0<r<1$ the set of real $t$ such that $re^{it} \in E$ has measure zero.
\end{thm}

Note that if $\phi$ and $\theta$    are inner then $ \phi  \circ \theta$ is also inner (this is not obvious). Here we are considering simply $b \circ \theta$ where $b$ is the inner function
with $b(z)=\dfrac{z-\alpha}{1-\overline\alpha z}$.\\

Frostman gave a stronger version of his theorem, expressed by saying that
the exceptional set has  
logarithmic capacity zero; however, it is beyond the scope of this work.

\begin{thm}
The {\em Crofoot transform}, defined for $\alpha \in \DD$ by
\[
J_\alpha f = \frac{(1-|\alpha|^2)^{1/2}}{1-\overline\alpha \theta} f \qquad (f \in K_\theta),
\]
is a unitary mapping from $K_\theta$ onto $K_{b \circ \theta}$ for each inner function $\theta$.
\end{thm}
In combination with Frostman's theorem, this can be used to construct orthonormal bases for any model space $K_\theta$.\\

\subsection{Truncated Toeplitz and Hankel operators}

Truncated Toeplitz operators were introduced by Sarason~\cite{sarason07}, and have
received much attention since then.
The idea here is to put finite Toeplitz matrices of the form
\beq\label{eq:fintoep}
\begin{pmatrix}
a_0 & a_{-1} & \ldots & a_{-n} \\
a_{1} & a_0 & \ldots & a_{-n+1}\\
\ldots & \ldots & \ldots & \ldots\\
a_{n} & a_{n-1} & \ldots & a_0
\end{pmatrix}
\eeq
into a more general context. One may also consider finite Hankel matrices of the form
\beq\label{eq:finhank}
\begin{pmatrix}
a_{-1} & a_{-2} & \ldots & a_{-n-1} \\
a_{-2} & a_{-3} & \ldots & a_{-n-2}\\
\ldots & \ldots & \ldots & \ldots\\
a_{-n-1} & a_{-n-2} & \ldots & a_{-2n-1}
\end{pmatrix}.
\eeq
Take $\theta$ inner, and $\phi \in L^\infty(\TT)$; then the truncated
Toeplitz operator 
$
A^\theta_\phi: K_\theta \to K_\theta
$
is defined by 
\[
A^\theta_\phi f = P_{K_\theta} (\phi \cdot f) \qquad (f \in K_\theta),
\]
where $P: L^2(\TT) \to K_\theta$ is the orthogonal projection.\\

 The motivating example involves the
choice $\theta(z)=z^{n+1}$, and the orthonormal basis $\{1,z,z^2,\ldots,z^n\}$ of $K_\theta$, when
the matrix of $A_\phi^\theta$ has the form (\ref{eq:fintoep}), with $(a_n)_{n \in \ZZ}$ the 
Fourier coefficients of $\phi$.\\

Similarly for truncated Hankel operators. The operator $B_\phi^\theta: K_\theta \to \overline{z K_\theta}$ is
defined by 
\[
B^\theta_\phi f = P_{\overline{zK_\theta}} (\phi \cdot f) \qquad (f \in K_\theta).
\]
Now, if $\theta(z)=z^{n+1}$, then $\overline{zK_\theta}$ has basis $\{\overline z,\ldots,\overline{z}^{n+1}\}$,
and with these bases the operator $B^\theta_\phi$ has a truncated Hankel matrix
(\ref{eq:finhank}).\\

\section{Restricted shifts}

\subsection{Basic ideas}

We recall that the invariant subspaces of the backwards shift $S^*$  have the form $K_\theta$.
We now define $S_\theta: K_\theta \to K_\theta$ by 
\[
S_\theta = P_{K_\theta} S_{|K_\theta} = (S^*_{|K_\theta})^*.
\]
This is the truncated Toeplitz operator with symbol $z$, and if we take $\theta(z)=z^{n+1}$ it maps
as follows: $1 \mapsto z,  z \mapsto z^2, \ldots, z^{n-1} \mapsto z^n, z^n \mapsto 0$, so that
its matrix is given by
\[
\begin{pmatrix}
0 & 0 & 0 & \ldots & 0 & 0 \\
1 & 0 & 0 & \ldots & 0 & 0\\
0 & 1 & 0 & \ldots & 0 & 0\\
\ldots & \ldots & \ldots & \ldots & \ldots & \ldots \\
0 & 0 & 0 & \ldots & 0 & 0 \\
0 & 0 & 0 & \ldots & 1 & 0
\end{pmatrix}.
\]

The restricted shift has a part in the Sz.-Nagy--Foias functional model~\cite{SNFBK}: if $T$ is a  
contraction on a Hilbert space $H$ such that
$\|(T^*)^n x \| \to 0$ for all $x \in H$ and $\rank(I-T^*T)=\rank(I-TT^*)=1$,
then there is an inner function $\theta$ such that $T$ is unitarily equivalent to $S_\theta$.

\begin{prop}
The invariant subspaces for the restricted shift $S_\theta$ are ``shifted'' model spaces of the form $K_\theta \cap \phi H^2=\phi K_{\theta/\phi}$, where $\phi$ is an inner function dividing $\theta$ in $H^\infty(\mathbb{D})$.
\end{prop}

\beginpf
The invariant subspaces for its adjoint, $S^*_{|K_\theta}$ are clearly of the form $K_\phi$, where
$\phi$ divides $\theta$ in $H^\infty(\mathbb{D})$. Their orthogonal complements are the invariant subspaces for $S_\theta$, and have
the required form.
\endpf

It is easy to see that $\rank S_\theta < \infty$ if and only if $\theta$ is a finite Blaschke product. 
We now define the {\em spectrum\/} of an inner function $\theta$ by
\[
\sigma(\theta) = \{w \in \overline \DD: \liminf_{z \to w} |\theta(w)|=0\}.
\]
For a Blaschke product $B$, the set $\sigma(B)$ is the closure of the
zero set of $B$ in $\overline{\DD}$.
It can then be shown that in general
$\sigma(S_\theta)=\sigma(\theta)$ (see \cite[Lec. VIII]{helson}).

\subsection{Unitary perturbations and dilations} 
We shall now suppose that $\theta(0)=0$: this simplifies some
of the formulae, but is not a serious restriction.
D.N. Clark \cite{clark} initiated a very fruitful study of unitary perturbations of restricted shifts.
In particular, he showed that the set of rank-1 perturbations of $S_\theta$ that are unitary
can be parametrised as $\{U_\alpha: \alpha \in \TT \}$,
where
\[
U_\alpha f = S_\theta f + \alpha \langle f, S^*\theta \rangle 1, \qquad (f \in K_\theta),
\]
noting that the constant function $1$ lies in $K_\theta$ because $\theta(0)=\langle\theta,1\rangle=0$.\\

If we consider the case $\theta(z)=z^{n+1}$, as above, we find that the matrix of $U_\alpha$ is now
\[
\begin{pmatrix}
0 & 0 & 0 & \ldots & 0 & \alpha \\
1 & 0 & 0 & \ldots & 0 & 0\\
0 & 1 & 0 & \ldots & 0 & 0\\
\ldots & \ldots & \ldots & \ldots & \ldots & \ldots \\
0 & 0 & 0 & \ldots & 0 & 0 \\
0 & 0 & 0 & \ldots & 1 & 0
\end{pmatrix},
\]
so that $1,z, \ldots,z^{n-2},z^{n-1}$ are mapped, respectively, to
$z,z^2,\ldots,z^{n-1},\alpha$.\\

The spectral measure of $U_\alpha$ is called a {\em Clark measure}, and there are various applications.
See, for example, the book \cite{CMR}.\\

 For  an operator $T$ on a  Hilbert space $H$, we consider the question of finding a unitary operator $U$ on a space containing $H$, such that
its restriction to $H$ is $T$. In matrix terms we may write
\[
U= \begin{pmatrix}
T & * \\
* & *
\end{pmatrix}.
\]
If $U$ is defined on $H \oplus \CC$, then we call it a 1-dilation. This is not the same as the
standard Sz.-Nagy--Foias dilation as in~\cite{SNFBK}.
In the context of restricted shifts and unitary dilations, there is a connection here with a classical result in geometry, which we now develop.

\subsection{Numerical ranges}

\begin{figure}
\begin{center}
\vskip-5.5cm
 \includegraphics[scale=0.5,angle=0]{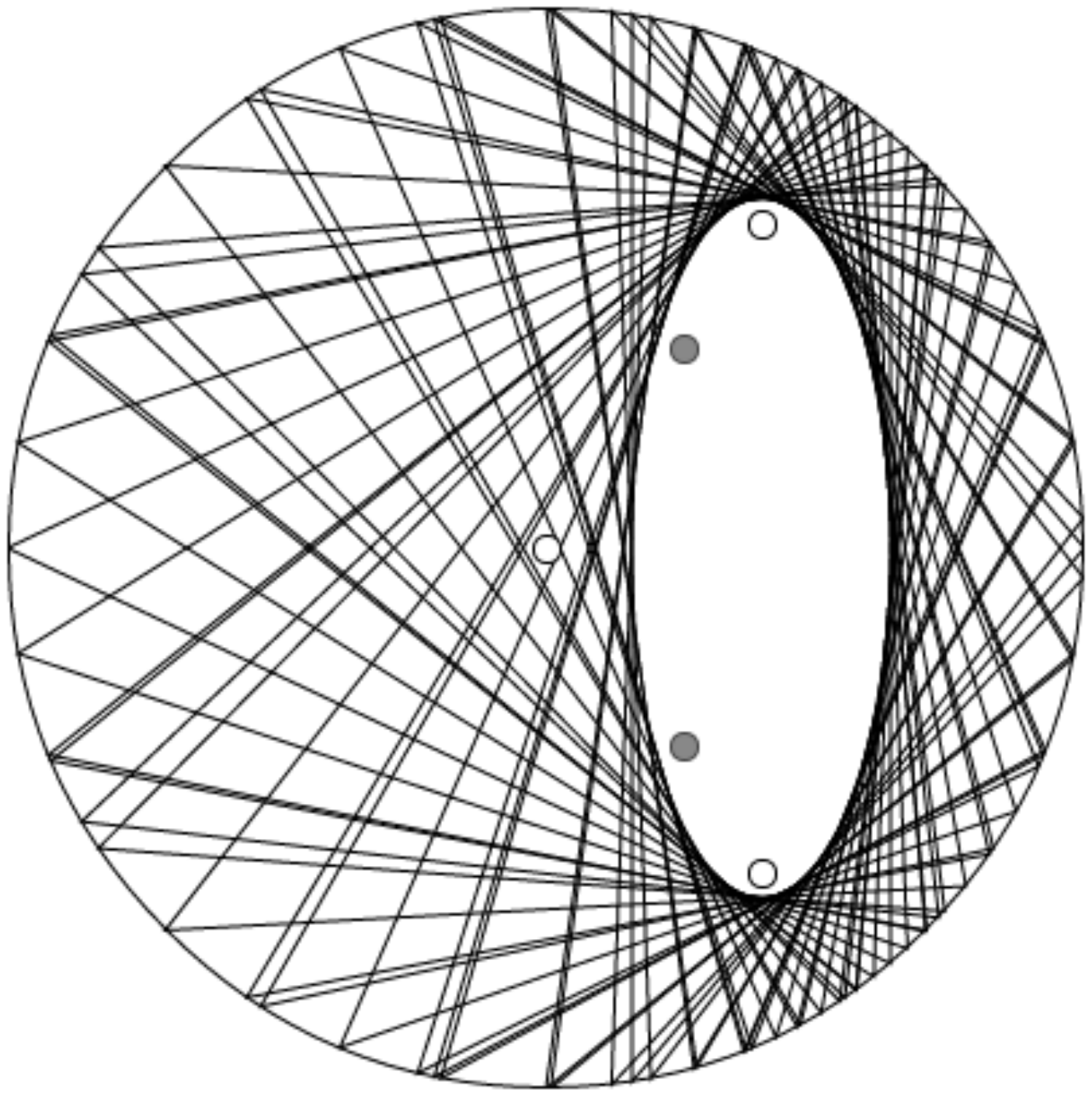}
\vskip-3.5cm
\caption{Poncelet ellipse with triangles}
\label{fig:1}
\end{center}
\end{figure}

For an integer $n \ge 3$, a closed subset $A$ of $\DD$ has the {\em $n$-Poncelet property}, if 
whenever there exists an $n$-gon $P$ such that $P$ circumscribes $A$ and 
has its vertices on
$\mathbb{T}$, then every point on the unit circle is a vertex of such an $n$-gon.
This was originally studied in the context of an ellipse, as in
Figure~\ref{fig:1}. (The figures were produced by an applet written by A.~Shaffer.)
Associated with the ellipse is a Blaschke product, as we shall explain: its zeroes are denoted by
light circles and the zeroes of its derivative by dark circles.

We shall also be considering a generalization of this, namely, an infinite Poncelet property.

Let us suppose first that $\theta$ is a finite Blaschke product, and hence $K_\theta$ is finite-dimensional.
Recall that the {\em numerical range\/} of an operator $T$ on a Hilbert space $H$ is defined by
\[
W(T)= \{ \langle Tx,x \rangle \}: x \in H, \|x\|=1 \},
\]
and, according to  the Toeplitz--Hausdorff theorem, is a convex subset of the plane. If $T$ has finite rank, then $W(T)$ is also compact.

\begin{thm} For the restricted shift $S_\theta$ on a finite-dimensional model space $K_\theta$ we have
\[
W(S_\theta) = \bigcap_{\alpha \in \TT} W(U^\theta_\alpha),
\]
where the $U^\theta_\alpha$ are the rank-1 Clark perturbations of $S_{z\theta}$, which are equivalent to unitary
1-dilations of $S_\theta$. 
\end{thm}

For versions of this results and further developments, see \cite{GW98,GW03,GR08,DGV10}.\\

Note that 
\[
\sigma(U^\theta_\alpha)=\{z \in \TT: z\theta(z)=\alpha\},
\]
an $n+1$-point set if the degree of $\theta$ is $n$. Moreover, $W(U^\theta_\alpha)$ is the convex hull
of $\sigma(U^\theta_\alpha)$, namely, a polygon.
If $\deg \theta=2$, then it is known that $W(S_\theta)$ is an ellipse, 
with foci at the eigenvalues of $S_\theta$. Therefore,
this ellipse has
foci at the zeroes of $\theta$, and it is here expressed as an intersection of triangles.\\

\begin{figure}
\begin{center}
\vskip-5.5cm
 \includegraphics[scale=0.5,angle=0]{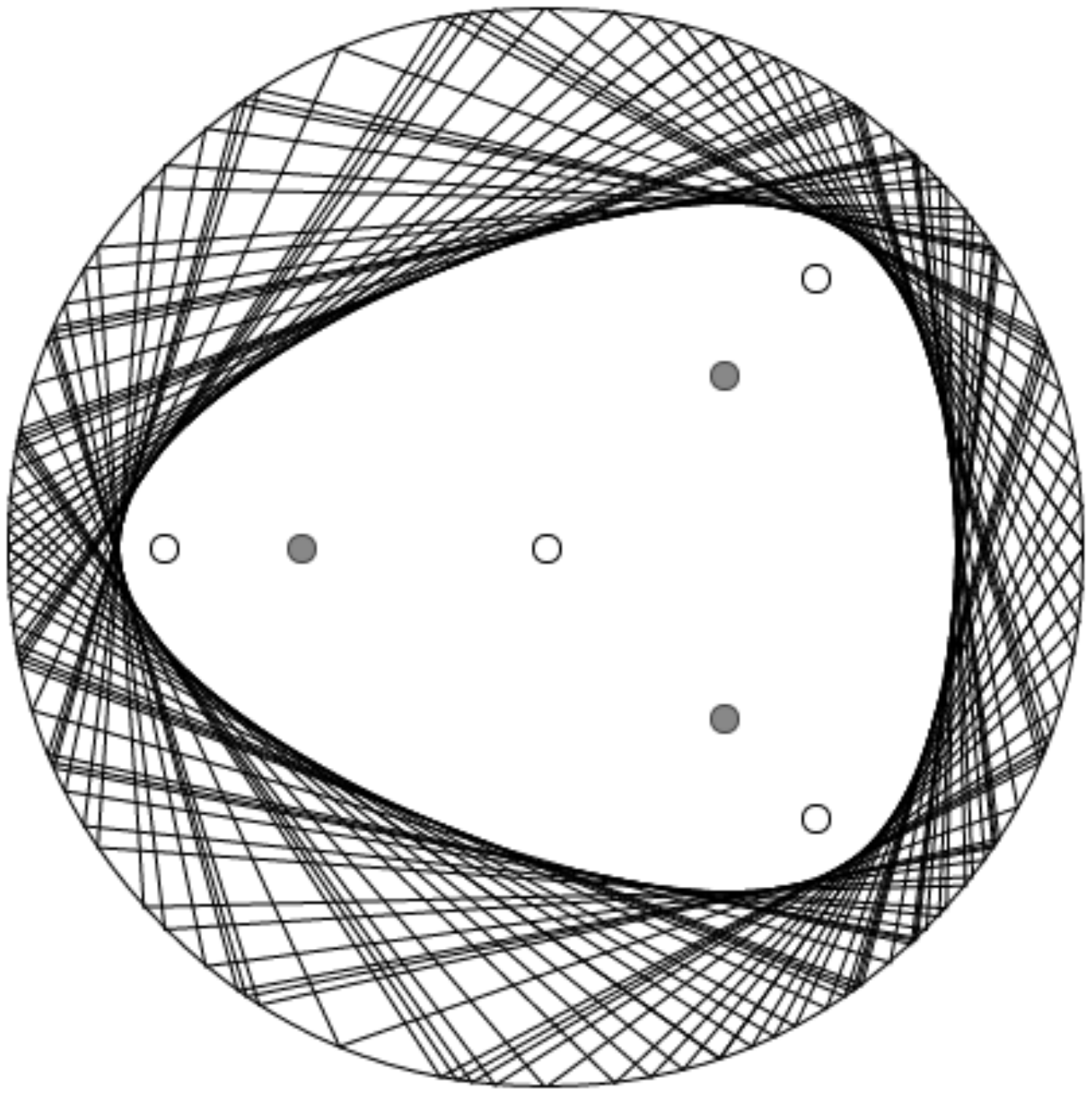}
\vskip-3.5cm
\caption{Symmetrical Poncelet curve with quadrilaterals}
\label{fig:2}
\end{center}
\end{figure}

\begin{figure}
\begin{center}
\vskip-5.5cm
 \includegraphics[scale=0.5,angle=0]{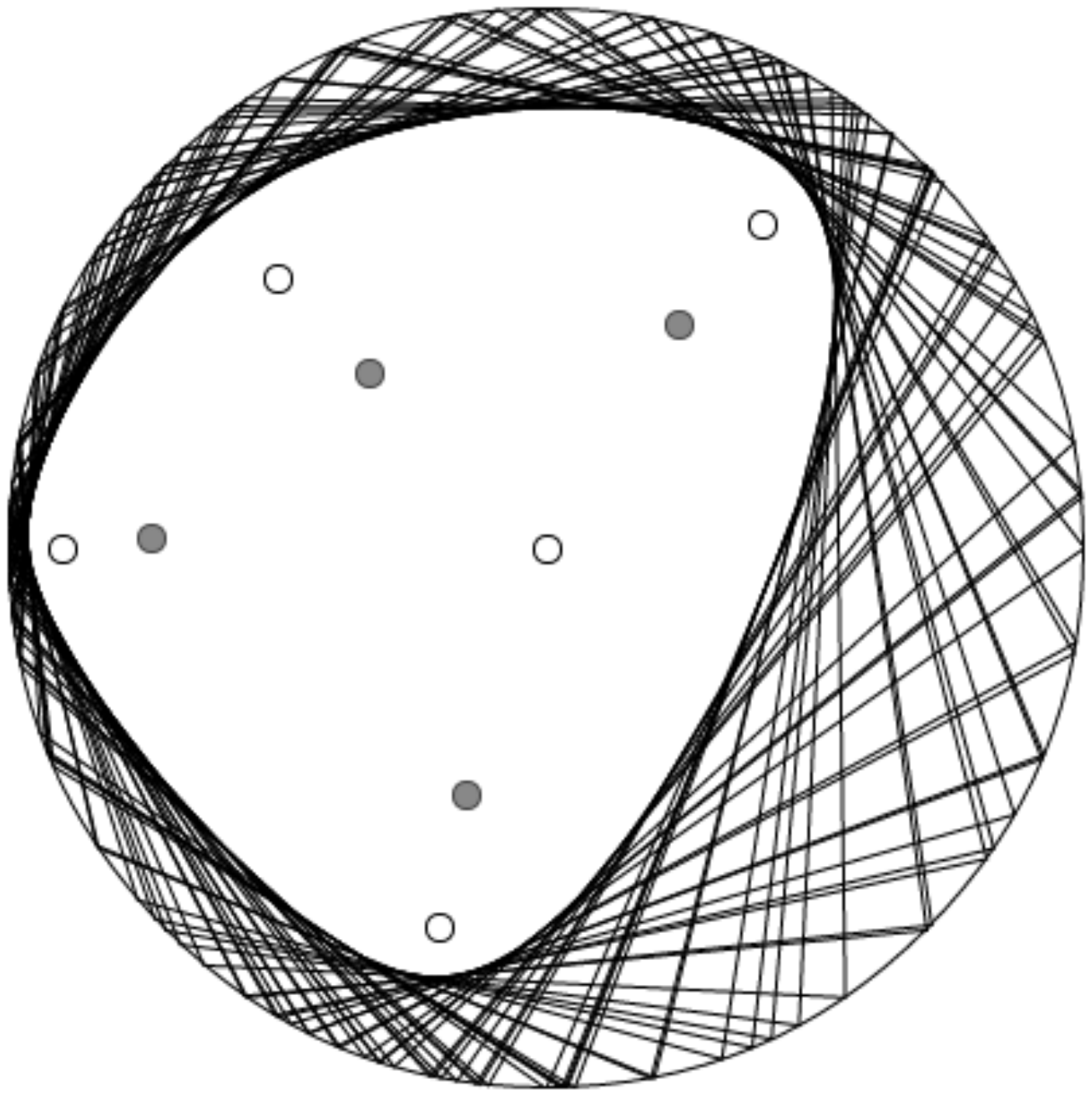}
\vskip-3.5cm
\caption{Asymmetrical Poncelet curve with pentagons}
\label{fig:3}
\end{center}
\end{figure}

Figures~\ref{fig:2} and \ref{fig:3} show similar examples with $n=3$ (quadrilaterals) and $n=4$ (pentagons).

The following more general  result was proved in \cite{CGP09}. Note that numerical ranges no longer need to be closed, so the 
formulation is slightly different.

\begin{thm} Let $\theta$ be an inner function. Then
\[
\overline{W(S_\theta)}= \bigcap_{\alpha \in \TT} \overline{W(U^\theta_\alpha)},
\]
where the $U^\theta_\alpha$ are the unitary 1-dilations of $S_\theta$ (or, equivalently, the rank-1 Clark perturbations of $S_{z\theta}$).
\end{thm}

In general we may regard the numerical ranges of the $U^\theta_\alpha$ as convex polygons with infinitely-many sides.
Some vectorial generalizations of these results (involving more general
contractions) are given in \cite{BGT11,BT14}.

We may now ask  how many polygons are needed to determine $\theta$ uniquely.
Note that the vertices of a polygon are solutions to $z\theta(z)=\alpha$, so we are motivated to consider boundary interpolation by inner functions. 

\subsection{Interpolation questions}

For finite Blaschke products we have the following theorem in \cite{CGP11}
about identifying two sets of $n$ points. Note that the two sets $\{z_1,\ldots,z_n\} $ and
$\{w_1,\ldots, w_n\}$ in the  theorem are necessarily interlaced; that is, each $z_j$ lies
between two successive $w_k$ and vice-versa.
\begin{thm}
For a finite Blaschke products $\theta$, $\phi$ of degree $n$, suppose that there are distinct points $z_1,\ldots,z_n $ and
$w_1,\ldots, w_n$ in $\TT$ such that
\begin{eqnarray*}
&\theta(z_1)=\ldots=\theta(z_n), \qquad &\theta(w_1)=\ldots=\theta(w_n),\\
\noalign{\hbox{and}}
&\phi(z_1)=\ldots=\phi(z_n), \qquad &\phi(w_1)=\ldots=\phi(w_n).
\end{eqnarray*}
Then $\phi = \lambda\dfrac{\theta-a}{1-\overline a \theta}$
for some $\lambda \in \TT$ and $a \in \DD$.
\end{thm}

We say that $\phi$ is a {\em Frostman shift\/} of $\theta$.\\

\begin{figure}[htbp]
  \begin{center}
    \leavevmode
\includegraphics{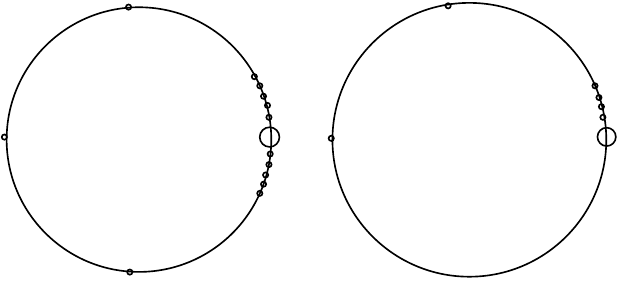}
    \caption{\emph{Singularities of type $2$ (L) and type $1$ (R)}}
    \label{fig:page8}
  \end{center}
\end{figure}

Suppose now that $\theta$ is inner with just one singularity on $\TT$; this is, it extends analytically across $\TT$ except at one point, which we
shall take to be $z=1$. For some such $\theta$, but not all, there will be a sequence $(t_n)_{n \in \ZZ}$ in $\TT$ (necessarily isolated since $\theta$ has an analytic extension),
accumulating on both sides of the point $1$, such that $\theta(t_n)=1$ for each $n$. 
This is called a singularity of Type 2 in \cite{CGP12}: see Figure~\ref{fig:page8}.\\

We consider how to determine $\theta$ from this data.

We transform to the upper half-plane $\CC^+$, using the M\"obius mapping
\[
\psi(z)=i\frac{1+z}{1-z}, \qquad \hbox{with} \quad \psi(1)=\infty.
\]
Now consider $F:= \psi \circ \theta \circ \psi^{-1}$. Then $F$ is meromorphic on $\CC$ with real poles $(b_n)_{n \in \ZZ}$ accumulating at $\pm\infty$. It maps
$\CC^+$ to $\CC^+$ and $\CC^-$ to $\CC^-$. Such functions are called {\em strongly real}. Without loss of generality we may assume that $0$ is neither a pole nor a zero of $F$,
in which case we have the following theorem, given in \cite{levin} as the Hermite--Biehler theorem, but attributed to Krein.
  
\begin{thm}
For $F$ strongly real with poles $(b_n)$ tending to $\pm\infty$,
the zeroes $(a_n)$ and poles $(b_n)$ are interlaced in the sense that we may write $b_n<a_n<b_{n+1}$ for each $n$, and then
\beq\label{eq:firstF}
F(z)=c \prod_{n \in \ZZ} \frac{1-z/a_n}{1-z/b_n},
\eeq
where $c>0$ unless $a_nb_n < 0$, in which case $c<0$. There is such a function for each sequence $(a_n)$ interlaced with the $(b_n)$.
\end{thm}

Our conclusion is that, given one limit point on $\TT$, approached from both sides by solutions
to $\theta(z)=1$, the set $\theta^{-1}(1)$ does not determine $\theta$, whereas
the sets $\theta^{-1}(1)$ and $\theta^{-1}(-1)$ together tell us what $\theta$ is, to within
composition by a M\"obius transformation fixing $\pm 1$. \\

In \cite{CGP11} the case of finitely-many singularities is discussed, including cases then
some singular points are approached on one side only. Curiously, there is a non-uniqueness
case in the Hermite--Biehler expression, apparently missed by Krein. For suppose
that $a_n \to 1$ as $n \to -\infty$ and $a_n \to \infty$ as $n \to \infty$. Then, with interlaced $(b_n)$
there is one solution, namely (\ref{eq:firstF}), but there is also another possibility, namely
\[
F(z)=c(z-1) \prod_{n \in \ZZ} \frac{1-z/a_n}{1-z/b_n}
\]
and these are the only possibilities. \\

On the circle, the corresponding $\theta$  has   a singularity of Type 1  in the terminology of \cite{CGP12}: see Figure~\ref{fig:page8}.
Thus there are two one-parameter families of
inner functions $\theta$ for such a choice of $\theta^{-1}(1)$ and $\theta^{-1}(-1)$.
A third set, e.g. $\theta^{-1}(i)$, enables one to distinguish between them. Thus
one sees that, in a fairly general situation, if $W(S_\theta)=W(S_\phi)$, then $\theta$ is a Frostman
shift of $\phi$ and so the restricted shifts are unitarily equivalent.\\

Some (necessarily less explicit) extensions of these ideas have been given by
Bercovici and Timotin \cite[Cor.6.3]{BT12}, in the case where the set
of singularities of the inner function $\theta$ is of measure zero.

\end{document}